\documentclass[12pt]{amsart}
\usepackage{amsmath,amssymb, euscript, graphicx}

\theoremstyle{definition}

\newcommand \smm {\displaystyle \sideset {}{'}\sum }
\newcommand \Cal {\mathcal }

\newcommand \dl {\delta}

\newcommand \ep {\epsilon}

\newcommand \al {\alpha}

\newcommand \bt {\beta}

\newcommand \ba {\bold a}
\newcommand \bb {\bold b}

\newcommand \bu {\bold u}

\newcommand \bw {\bold w}

\def \qed{\ifhmode\unskip\nobreak\fi\ifmmode\ifinner\else\hskip5pt\fi\fi
 \hbox{\hskip5pt\vrule width4pt height6pt depth1.5pt\hskip1pt}}
\begin{document}
\centerline {\bf GOOD MODULATING SEQUENCES} 
\centerline {\bf FOR THE ERGODIC HILBERT TRANSFORM} 
\bigskip
\bigskip

\centerline { AZER AKHMEDOV and DO\u GAN \c C\"OMEZ
\footnote {* {\it 2000 Mathematics Subject Classification.} Primary 28D05; Secondary 42B20, 40A30. \newline
}} 

\bigskip

\noindent{\bf ABSTRACT.} {\Small This article investigates classes of
bounded sequences of complex numbers that are universally 
good for the ergodic Hilbert transform in $L_p$-spaces, $2\leq p \leq \infty.$
The class of bounded Besicovitch sequences satisfying a rate condition
is among such sequence classes.}  

\bigskip

\noindent {\bf 1. Introduction.}  Let $(X,\Sigma ,\mu)$ be a 
probability space and $T:X \to X$ be an invertible measure 
preserving transformation.  For any $f\in L_p, $ the ergodic
Hilbert transform (eHt) of $f$ is defined as 
$$H f(x) := \lim_n \sideset {}{'}\sum_{k=-n}^{n} \frac{T^k f(x)}{k}, $$
if the limit exists, 
where $\smm_{k=-n}^{n} $ means summation without the term $k=0. $
It is well-known that the eHt exists a.e. 
for $f\in L_1 $ [C, P].  This result has also been extended to 
various other settings [A, CP, J, S, \c C$_1$].  
Given a sequence $\ba =\{a_k\} $ of complex
numbers, we will define the modulated ergodic Hilbert
transform of $f\in L_p $ (modulated by $\ba $) as
$H_{\ba} f(x) := \displaystyle \lim_n \smm_{k=-n}^{n} \frac{ a_k T^k f(x)}{k}. $  
If $(X,\Sigma ,\mu, T)$ is a dynamical system, a sequence $\ba $ is
called {\it good for the ergodic Hilbert transform in $L_p(X)$} if
the modulated ergodic Hilbert transform exists $\mu$-a.e. for every
$f\in L_p (X).$  
Let $\Cal{T} $ be a class of measure preserving dynamical systems. 
We will say that the sequence $\ba $ is {\it universally good for the 
ergodic Hilbert transform in $L_p$ for the class $\Cal{T} $} 
if $\ba $ is good for the ergodic Hilbert transform in $L_p$ of every 
dynamical system in $\Cal{T} . $ 
In case that $\ba $ is good for the ergodic Hilbert 
transform in $L_p$ of {\it every} dynamical system,
we will say that it is {\it universally good for the ergodic Hilbert
transform in $L_p$.}

In the article [\c C$_2$] the second author investigated 
some classes of sequences that are universally good for the eHt.  Such sequence classes 
are rather large; for instance, symmetric sequences of bounded variation 
and sequences of Fourier coefficients of functions in 
$L_p [0, 2\pi] ,\ 1<p<\infty , $ are universally good for the eHt.  
Recently, in [LT] a Wiener-Wintner type theorem for the ergodic Hilbert transform 
was proved; a remarkable result which eluded mathematicians for two decades.
A direct consequence of this result is that the sequences of the form
$\{\lambda^k \}, \ \text{with}\ |\lambda |=1, $ are universally good for the 
eHt in $L_p, \ 1<p<\infty . $
The techniques utilized in [\c C$_2$] fell short of proving that some sequence 
classes, known to be universally good for the ergodic averages, are universally 
good for the eHt.  There, besides indicating that not all
bounded Besicovitch sequences are good modulating for the eHt, 
it was proved that a proper subset of the set of bounded Besicovitch 
sequences is good for the eHt in $L_1(X) . $
In this article, having the Wiener-Wintner theorem for the ergodic Hilbert 
transform, we will prove that those sequences 
are universally good for the eHt in $L_2 . $  Since we are in a probability 
space setting, these results also hold for $L_p$-functions for $2\leq p < \infty . $ 
We also obtain other classes of sequences universally good for the eHt.  
Throughout this article, unless stated otherwise, we will assume $0<\bt <1 $  
and $ 1<\al \leq 2 . $  Also, $\Delta $ will denote the unit circle in complex plane, and $C$ will always denote a constant, which may not be
the same at each occurrence.
\bigskip

\noindent{\bf 2. Bounded universally good sequences for the eHt.}
In this section we will show that some fairly large classes of bounded 
complex sequences that satisfy a rate condition are universally good for the eHt.  
Let $\ba $ be a sequence such that
$$ (\ast ) \ \ \sum_{k=-n}^{n} |a_k | =O(n^{\bt}), \ \ (n\geq 1) ,$$
In [\c C$_2$] it is shown that if $\ba $ is a bounded sequence good for the 
ergodic theorem in $L_{\infty} $ and satisfy the condition $(\ast )$, then it 
is universally good for the eHt in $L_1 . $
The class of sequences satisfying $(\ast )$ include the sequences of Fourier coefficients
of functions belonging to the function spaces $L_p [0, 2\pi],\ 1<p<\infty , $
$L_{\al} [0, 2\pi] $ (the ${\al}-$Lipschitz functions in
$L_1 [0, 2\pi] $), and $BV_1 [0, 2\pi] $ (the functions of bounded variation
in $L_1 [0, 2\pi]$).  The condition $(\ast )$ is naturally satisfied by the 
sequences belonging to these classes; however, the same assertions
made there are also valid if one considers sequences satisfying a
weaker condition.  For, define
$$M_{\al}=\{ \ba : \sum_{k=-n}^{n} |a_k | =O(\frac{n^{\al -1}}{\log^{\al}n}) \} . $$
If $n$ is large enough and $\al > \bt + 1$, then $n^{\bt} \leq \frac{n^{\al -1}}{\log^{\al}n} $; hence,
any sequence satisfying the condition $(\ast )$ belongs to $M_{\al} $,
if $\al > \bt + 1.$  

In [CL] it is shown that, among several other results related to the 
{\it one sided} ergodic Hilbert transform, if $\ba =\{a_k\}_{k\in \Bbb Z} $ 
is a sequence of complex numbers satisfying
$$ (\ast \ast ) \ \ \sup_{n\geq 1} \max_{|z |=1} \frac{1}{n^{1-\bt}}
|\sum_{k=1}^{n} a_{k} z^{k} | = C_{\ba} <\infty ,$$
then it is universally good for the eHt in $L_1 . $  When the (two-sided)
ergodic Hilbert transform is concerned, it turns out that one can consider
a larger class of sequences.
Let $A_{\al}$ denote the set of sequences $\ba =\{a_k\}_{k\in \Bbb Z} $ 
of complex numbers satisfying
$$ (1) \ \ \sup_{n\geq 1} \max_{|z |=1} \frac{\log^{\al}n}{n^{\al -1 }}
|\sum_{k=-n}^{n} a_{k} z^{k} | = C_{\ba} <\infty .  $$
Since $\frac{n^{\al - 1}}{\log^{\al}n} \geq n^{1-\bt} $ when $n$ is sufficiently large and $\al > 2-\bt , $
all classes $A_{\al}$ contain sequences satisfying the condition $(\ast \ast )$.

\noindent{\bf Remarks.} 1.  All the sequence classes mentioned above do not 
contain constant sequences; on the other hand, constant sequences are trivially 
universally good for the eHt [C,P]. 

2. $M_{\al} \subset A_{\al} . $  In fact, $A_{2} $ contains all $M_{\al} $ for all 
$1<\al \leq 2 . $

3. If $\ba \in M_{\al}, $ then for all $f\in L_{\infty}, $
$ \lim_{n\to \infty} \frac1n \sum_{k=0}^{n-1} a_k T^k f =0 $ a.e. 

4. If $\al > 1+\bt , $ then $\frac{n^{\al -1 }}{\log^{\al}n} \geq n^{\bt}. $
Hence it follows that any sequence satisfying the condition $(\ast )$ belongs to
$A_{\al}. $ In particular, $A_2 $ contains all such sequences.

5. There are $A_{\al}$ sequences that do not belong to any $M_{\al} . $
For instance, let $\ba =\{a_n \} $ be the special case of the Hardy-Littlewood
sequence given by $ a_n=e^{in \log |n|}. $  Clearly, 
$ \ba \notin M_{\al} $ since $\displaystyle \sum_{k=-n}^{n} |a_k |= O(n). $ 
However, for any $|z|=1, $ it follows that
$ | \displaystyle \sum_{k=-n}^{n} a_k z^k | = O(\sqrt{n}) $ (see [Z, p: 199]);
hence $ \ba \in A_{3/2} . $ 
\medskip

A sequence satisfying $(\ast )$(hence in $M_{\al} $)
need not be bounded.  For example, let $a_k=j $ if $k=\mp 2^j , $ and $a_k=0 $
for otherwise.  Then $\ba $ is an unbounded sequence and satisfies $(\ast )$ for any $\bt\in (0,1)$.  Having this note, however, all the sequences considered
throughout the rest of this article will be bounded.

\medskip

Although sequences satisfying the condition $(\ast )$ are included in $A_{2}$, and $M_{\al} \subset A_{\al} , $ for some values of $\al $ we also have the reverse inclusion.

\noindent{\bf Proposition 2.1}  {\it Let $\al^{\prime} +1/2 < \al \leq 2 , $ then
$A_{\al^{\prime}} \subset M_{\al} . $}

\noindent{\bf Proof.}  By H\"older's inequality,
$$\sum_{k=-n}^{n} |a_k | \leq (2n)^{1/2} (\sum_{k=-n}^{n} |a_k |^2)^{1/2}
=(2n)^{1/2} [\int_{\Delta} |\sum_{k=-n}^{n} \lambda^k a_k |^2 d\lambda ]^{1/2}. $$
Since $ \ba \in A_{\al^{\prime}} , $ we have
$\sum_{k=-n}^{n} |a_k | \leq C \frac{n^{\al^{\prime} -1/2}}{\log^{\al^{\prime}}n } $ for 
some constant $C. $ Hence, since $ \al^{\prime} +1/2 < \al \leq 2 , $ the assertion 
follows. \qed
\medskip

\noindent{\bf Remark.}  It follows from Proposition 2.1 that, if
$\ba \in A_{\al}, \ 1<\al < 3/2, $ then $3/2 <1/2+ \al < 2 , $ and hence,
for all $f\in L_{\infty}, $
$ \displaystyle \lim_{N\to \infty} \frac1N \sum_{k=0}^{N-1} a_k T^k f =0 $ a.e.
\medskip

Next, we will prove that $M_{\al}$ sequences are universally good for the eHt.
The proof is essentially the same as the proof of Theorem 2.2 in [\c C$_2$]; hence,
we will sketch it here for completeness.

\noindent{\bf Theorem 2.2} {\it Let $\ba=\{a_k \} \in M_{\al} . $ Then we have
the weak (1,1) maximal inequality for $H_{\ba} f $: for any 
$f\in L_1, $ and for any $\lambda >0 , $ there is a constant $C$ such that
$$ \mu(\{ x: \sup_{n\geq 1 } \vert
\smm_{k=-n}^{n}\frac{a_k T^kf(x)}{k} \vert > \lambda \})
\leq \frac{C}{\lambda } \Vert f \Vert_1 . $$ 
Furthermore, $\ba $ is universally good for the eHt in $L_1 . $} 

\noindent{\bf Proof.}  By Abel's summation by parts formula,
$$ \smm_{k=-n}^{n} \frac{a_k T^kf}{k} =
 \sum_{k=1}^{n-1} \frac{S_k - S_{-k}}{k(k+1)}
+ \frac1n (S_{n}- S_{-n} ) , \ \text{where}\
S_{\mp j}=\sum_{i=1}^{j} a_{\mp i} T^{\mp i}f . $$
If $E=\{ x: \displaystyle \sup_{n\geq 1 } \vert
\smm_{j=-n}^{n}\frac{ a_j T^jf(x)}{j} \vert > \lambda \}, $
then $E\subset E_1 \cup E_2 \cup E_3 , $ where
$E_1 =  \{ x: \displaystyle \sup_n |\frac1n S_n(x)|> \frac{\lambda}{3} \} , $
$E_2= \{ x: \displaystyle \sup_n |\frac1n S_{-n} (x)| > \frac{\lambda}{3} \} , $ and
$E_3= \{ x: \displaystyle \sup_n|\sum_{j=1}^{n}
\frac{1}{j(j+1)} [S_j(x)-S_{-j}(x)] |
> \frac{\lambda}{3} \}. $
Since we always have a weak (1,1) maximal inequality for the
operators $\frac{1}{n} S_{\mp n} $ when $\ba $ is a bounded sequence,
$\mu (E_1) \leq \frac{C_1}{\lambda} \Vert f \Vert_1 $ and
$\mu (E_2) \leq \frac{C_2}{\lambda} \Vert f \Vert_1, $
for some constants $C_1$ and $C_2. $

If $f\in L_1 $ and $\ba \in M_{\al}, $ we have, for some constant $C, $
$$ \int |\sum_{1\leq k \leq n} \frac{S_k-S_{-k}}{k(k+1)} |
\leq  \Vert f \Vert_1 \sum_{1\leq k \leq n} \frac{1}{k^2 } \sum_{j=-k}^{k}
|a_j | \leq \Vert f \Vert_1
\sum_{1\leq k \leq n} \frac{C}{k^{3-\al} \log^{\al}k}
\leq C \Vert f \Vert_1,  $$
where $C= \displaystyle \sum_{1\leq k \leq \infty} \frac{C}{k^{3-\al} \log^{\al}k}. $
Since the sequence 
$\{h_n\}=\{ \displaystyle \sum_{1\leq k \leq n} \frac{1}{k^2 } \displaystyle \sum_{j=-k}^{k}
|a_j | T^j |f| \} \subset L_1 $ is monotone increasing with 
$\int h_n \leq C \Vert f \Vert_1, $ by the Monotone Convergence Theorem, 
$\int  h_n \uparrow \int h \leq C \Vert f \Vert_1$ where $h$ is the pointwise limit of the sequence $h_n$.   Hence, 
by Chebyshev's inequality, for any $\lambda>0 , $
$ \mu(E_3) \leq \frac{C}{\lambda } \Vert f \Vert_1 . $
Hence, the weak (1,1) maximal inequality for 
$H_{\ba} f $ follows.

If $f \in L_{\infty} $ and $m<n $ are positive integers, then
$$ |\sum_{m\leq k \leq n} \frac{S_k-S_{-k}}{k(k+1)} |
\leq \Vert f \Vert_{\infty} \sum_{k=m}^{n} \frac{1}{k^2 }
\sum_{j=-k}^{k} |a_j| \leq  \Vert f \Vert_{\infty}
 \sum_{k=m}^{n} \frac{C}{k^{3-\al} \log^{\al}k}, $$
which implies that the sequence
$\{ \displaystyle \sum_{1\leq k \leq n} \frac{1}{k(k+1)} (S_k-S_{-k})(x) \}$
is Cauchy a.e.; hence, it converges.  Since
$\displaystyle \lim_n \frac{1}{n} S_{\mp n}(x) $ also converges a.e. for
all $f \in L_{\infty} , \ \displaystyle \lim_n H_{\ba} f (x) $ exists a.e. for
all $f \in L_{\infty} . $  By the Banach Principle, this fact combined with
the weak (1,1) maximal inequality in the first part implies that 
$\ba $ is universally good for the eHt in $L_1 . $ \qed
\medskip

By Theorem 2.2 and Proposition 2.1, any $\ba \in A_{\al }$ is also universally 
good for the eHt if $1< \al \leq 3/2 . $  For $3/2 < \al \leq 2 $ we need
different arguments.  Indeed, the statement below provides an argument valid for
all $\ba \in A_{\al }, \ 1< \al \leq 2 . $
\medskip

\noindent{\bf Theorem 2.3}  {\it If
$\ba=\{a_k \} \in  A_{\al}, \ 1 <\al \leq 2 , $ is a 
sequence good for the ergodic averages,
then it is universally good for the eHt in $L_2 . $ }

\noindent{\bf Proof.} 
Since $\displaystyle \lim_n \frac{1}{n} S_{\mp n} $ exists a.e.
by assumption, using Abel's partial summation, in order to show that
$\displaystyle \lim_n \smm_{j=-n}^{n} \frac{ a_j T^j f(x)}{j} $ exists
a.e. for all $f \in L_2, $ all we need to show is that
$\displaystyle \lim_n  \sum_{j=1}^{n-1} \frac{S_j - S_{-j}}{j(j+1)} $
exists a.e., where $S_{\mp j}=\displaystyle \sum_{k=1}^{j} a_{\mp k} T^{\mp k} f(x). $

Observe that,
$ \Vert S_j -S_{-j} \Vert_2^2 =<S_j, S_j>+<S_{-j}, S_{-j}>-
<S_j, S_{-j}>-<S_{-j}, S_j >, $
where $<f,g>=\int f \overline{g} d\mu . $  Since
$$ <S_j, S_j> = \sum_{k,l=1}^{j} <a_k T^kf, a_l T^lf>
=\sum_{k,l=1}^{j} a_k \overline{a_l} <T^{k-l}f,f> , $$
by the spectral theorem for unitary operators, we have
$<T^{k-l}f,f>= \int_{\Delta} z^{k-l} d\mu_f(z), $ where $\Delta$
is the unit circle.  Therefore,
$ <S_j, S_j> = \int_{\Delta}[ \sum_{k,l=1}^{j}
(a_k z^k )(\overline{a_l z^l })] d\mu_f(z); $
and hence, it follows that
$$ \aligned \Vert S_j -S_{-j} \Vert_2^2 &=
\int_{\Delta} [\sum_{k,l=1}^{j} (a_k z^k \overline{a_l z^l}
-a_k z^k \overline{a_{-l} z^{-l}} - a_{-k} z^{-k} \overline{a_l z^l}
+a_{-k} z^{-k} \overline{a_{-l} z^{-l}})] d\mu_f(z) \\
& = \int_{\Delta}| \sum_{1\leq |k | \leq j}
a_k z^k |^2 d\mu_f(z) . \endaligned $$
Since  $ \ba \in A_{\al} , $ it satisfies (1); hence, we have
$| \displaystyle \sum_{-j}^{j} a_k z^{k} | \leq C_{\ba} \frac{j^{\al -1}}{\log^{\al} j} . $
Thus,
$ \Vert S_j -S_{-j} \Vert_2 \leq
C \frac{j^{\al -1}}{\log^{\al} j}  \Vert f \Vert_2, $
for some constant $C $ that depends on $\ba . $  Therefore, by H\"older's inequality, 
it follows that
$$ \int |\sum_{j=1}^{n-1} \frac{S_j - S_{-j}}{j(j+1)} |
\leq  \int \sum_{j=1}^{n-1} \frac{|S_j - S_{-j}|}{j(j+1)}
\leq C \Vert f \Vert_2 \sum_{j=1}^{n-1} \frac{1}{ j^{3-\al } \log^{\al}j} . $$
Now, by the Monotone Convergence Theorem
$$\int \displaystyle \lim_n |\sum_{j=1}^{n-1} \frac{S_j - S_{-j}}{j(j+1)} |\leq \lim_n \int \sum_{j=1}^{n-1} \frac{|S_j - S_{-j}|}{j(j+1)} 
\leq C \Vert f \Vert_2 \sum_{j=1}^{\infty} \frac{1}{ j^{3-\al } \log^{\al}j} <\infty; $$
hence, we deduce that 
$\displaystyle \lim_n \sum_{j=1}^{n-1} \frac{S_j - S_{-j}}{j(j+1)} $ exists a.e. \qed
\medskip

\noindent{\bf Remark.} It should be noted here that the arguments in
the theorem above are purely $L_2$ space arguments.  
\medskip

The class of Besicovitch sequences are known to be 
universally good for the ergodic averages [BL].   
In [\c C$_2$] it has been proved that the 
sequences of Fourier coefficients of functions in 
$L_p [0, 2\pi] ,\ 1<p \leq \infty , $ which are 
bounded Besicovitch sequences, are universally good for the eHt 
in $L_p, \ \ 1<p \leq \infty . $  There,
it was also observed that {\it not every} sequence $\ba \in B $ is
good modulating for the eHt.  On the other hand,
a smaller subclass $B_{\bt}$ of $B $ produces good modulating sequences
for the eHt, where
$$B_{\bt} =\{\ba\in l_{\infty} : \exists \ \bw \ \text{induced by
a trigonometric polynomial such that}\ \ba-\bw \ \text{satisfies (*)} \} . $$
The techniques used in [\c C$_2$], however, fell short of showing
that sequences in $B_{\al} $ are {\it universally} good for the eHt.
In this section, having the Wiener-Wintner theorem for eHt [LT], 
we will show that not only the sequence class $B_{\al} , $
but a subclass of bounded sequences that contains 
$B_{\al} $ provides sequences universally good for the eHt.

In [LT], is was shown that if $f\in L_p ,\ 1<p<\infty, $ then there is a set
$X_f \subset X $ of probability one such that for all $x\in X_f $
$$ \lim_n \smm_{-n}^{n} \frac{\lambda^k T^k f (x)}{k} \ \ \text{exists for all}\ \ 
| \lambda |=1 . $$
Let $\mathcal{W}$ denote the class of sequences induced
by bounded trigonometric polynomials, which are
finite linear combinations of sequences of the form
$ \{ \lambda^k \}, \ |\lambda |=1. $  Hence, it follows that
the Lacey-Terwilleger Theorem holds for sequences in $\mathcal{W} ; $
that is, if $\bw \in \mathcal{W}$ then it is 
universally good in $L_p, \ 1<p<\infty . $
\medskip

\noindent{\it Two-sided bounded Besicovitch sequences.}  First, 
we will consider two-sided
bounded Besicovitch sequences $\ba=\{a_k\}_{k\in \Bbb Z} \in B, $ 
which are defined as, 
given $\ep >0, $ there exists $\bw_{\ep} \in \mathcal{W} $ such that
$$(\dagger ) \ \ \limsup_n \frac1n \sum_{k=-n}^{n} |a_k -w_{\ep} (k)| < \ep. $$  
Now, we define 
$$\aligned  &MB_{\al}= \{ \ba \in l_{\infty} : \forall \ep>0 \ 
\exists \ \bw_{\ep} \in \mathcal{W} \ \text{such that}\
\limsup_n \frac{\log^{\al}n}{n^{\al-1}} \sum_{k=-n}^{n} |a_k - w_{\ep} (k) |
< \ep \} \ \text{and} \\
&AB_{\al}= \{ \ba \in l_{\infty}: \exists \ \bw \in \mathcal{W} \ \text{such that}\
\ba-\bw \in A_{\al} \}  .  \endaligned  $$
\medskip

\noindent{\bf Remarks.}  1. $MB_{\al} \subset AB_{\al} $ for all $\al . $

2.  If $\ba \in MB_{\al} , $ then $\ba $ is bounded
Besicovitch. 

3.  If $\ba \in B_{\bt}, $ then $a_k=w_k + b_k , $ where
$ \bw=\{w_k\} \in \mathcal{W}, $ and $\{b_k \} $ satisfies the condition $(\ast )$, hence $\ba-\bw \in M_{\al} . $ Therefore, 
$B_{\bt} \subset MB_{\al} \subset AB_{\al} $ for all 
$0<\bt <1 $ and $1< \al \leq 2 . $  

4.  Sequences induced by trigonometric polynomials belong to 
the sequence space $MB_{\al}; $ and hence, to $AB_{\al} . $ 
\medskip

By [LT] and the remarks above, any $\bw \in \mathcal{W} $ 
is universally good for the eHt in $L_2 . $
Since for $\ba\in AB_{\al} , $
$ a_k= a_k -w_{\ep}(k) + w_{\ep}(k) , $ where $\bw_k $ is the appropriate
trigonometric polynomial, all we need to prove is that
$\{ \smm_{-n}^{n} \frac{(a_k -w_{\ep}(k)) T^k f}{k} \}_n $ converges a.e.
In that case, using the same techniques as in Theorem 2.3,
we obtain 

\noindent{\bf Corollary 2.4}  {\it If $\ba \in AB_{\al} , $
then it is universally good for the eHt in $L_2 . $}
\medskip

\noindent{\bf Remark.} It follows from Corollary 2.4 that if 
$\ba=\{a_k \} \in B_{\bt}, \  0<\bt < 1, $
then it is universally good for the eHt in $L_2 . $ 
\medskip

\noindent{\it Symmetric (one-sided) bounded Besicovitch sequences.}
In this section we will consider symmetric bounded Besicovitch sequences,
namely, {\it ordinary} bounded Besicovitch sequences $\ba =\{a_k\} $ (with
$a_k=a_{-k} $) such that given $\ep >0, $ there exists $\bw_{\ep} \in \mathcal{W} $ 
satisfying
$$\limsup_n \frac1n \sum_{k=1}^{n} |a_k -w_{\ep} (k)| < \ep. $$

First, we make an observation.  Let $T:\Delta \to \Delta $ be an irrational 
rotation, say $Tz= \phi z $ for some $|\phi |=1,\ \ \phi \neq 1. $
Then for any $\lambda $ on the unit circle, and for any $f$ having
$\phi $ as eigenvalue,
$$ \smm_{k=-n}^{n} \frac{\lambda^{|k|} T^k f}{k} = f
\sum_{k=1}^{n} \frac{(\phi \lambda)^k - (\bar{\phi} \lambda)^k}{k}; $$
and hence, if $\lambda=\bar{\phi} , $ then this series is not convergent.  
Therefore, symmetric bounded Besicovitch sequences defined by 
trigonometric polynomials need not be good for irrational rotations, 
which is different form the two-sided case.

Consider sequences $\ba =\{a_k\}_{k\geq 0} $ such that
$\gamma_{\ba}(k) := \displaystyle \lim_n \frac1n \sum_{j=1}^{n} a_{j+k} \bar{a_j} $
exists for all $k\in \Bbb N. $  $\gamma_{\ba} $ is called the {\it correlation}
of $\ba , $ which is extended to negative integers by letting
$\gamma_{\ba} (-k)= \bar{\gamma_{\ba}(k)} . $  Sequences
$\{ \gamma_{\ba}(k) \} $ are positive definite; hence, by the Herglotz-Bochner
theorem there exists a unique Borel probability measure $\mu_{\ba} $
on the unit circle $\Delta $ such that
$$ \gamma_{\ba}(k) = \int_{\Delta} z^k d\mu_{\ba} (z), \ \ n\in \Bbb N. $$
The measure $\mu_{\ba} $ is called the {\it spectral measure of $\ba. $} 
Bounded Besicovitch sequences are known to have correlation; indeed,
bounded Besicovitch sequences are exactly those complex sequences such that:
(i) $ \mu_{\ba} $ is discrete, 
(ii) $ \Gamma(z):= \displaystyle \lim_n \frac1n \sum_{j=0}^{n} a_j \bar{z^j} $
exists for every $ z,$ and (iii) $ \mu_{\ba} (z) =|\Gamma(z) |^2 $
for all $ z\in \Delta $ [BL].  Furthermore, it is also known that 
$\Gamma(z)=0 $ for all but at most countably many $z\in \Delta $ [K]. 
The set $ \sigma (\ba) =\{ z\in \Delta : \Gamma(z) \neq 0 \} $
is called the {\it spectrum of $\ba . $}  Obviously, if $a_k=\lambda^k $
for some $\lambda \in \Delta, $ then $\sigma({\ba}) =\{\lambda \}. $ 

If $(X,\Sigma ,\mu, T)$
is an ergodic dynamical system, then $L_2(X)=\kappa \oplus \kappa^{\perp} , $
where $\kappa $ is the closed linear subspace spanned by the eigenfunctions
of $T $ (called the Kronecker factor of the system).  
Consequently, for a non-constant bounded Besicovitch sequence $\ba $
and a measure preserving system with $\{f\in L_2 : Tf=f\} \subset \kappa $
properly, if $ \lambda \in \sigma(\ba) \cap \sigma(T), $ then, as observed above,
$\displaystyle \lim_n \smm_{k=-n}^{n} \frac{\lambda^{|k|} T^k f}{k} $ need
not exist a.e.  These arguments also imply that, given any dynamical
system $(X, \Sigma, \mu, T) $ with $\{f\in L_2 : Tf=f\} \subset \kappa $ properly, 
there exists a bounded Besicovitch sequence $\ba $ such that
$\displaystyle \lim_n \smm_{k=-n}^{n} \frac{a_{|k|} T^k f}{k} $ fails to exist
a.e.  However, as the next result shows, symmetric sequences $\ba \in B_{\al} $ 
are still universally good in $L_2 $ in a restricted sense.
Let $\Im $ denote the class of weakly mixing measure preserving systems. 
If $T \in \Im , $ then it has continuous spectrum, and hence,
its Kronecker factor is simple; namely, $\kappa=\{f\in L_2 : Tf=f\}. $
Therefore, for a weakly mixing $T , $ 
$ \displaystyle \lim_n \smm_{k=-n}^{n} \frac{ \lambda^{|k|} T^k f}{k} $
exists a.e. for any $f\in \kappa . $  Again, 
using the same techniques as in Theorem 2.3,
we obtain
\medskip

\noindent{\bf Corollary 2.5}  {\it If $\ba \in AB_{\al} $ is a symmetric sequence, 
then it is universally good for the eHt in $L_2$ for the class $\Im . $}
\medskip 

The following theorem is analogous to Theorem 2.3, albeit with a restriction
on the class of transformations.
Let $ \mathcal{L} $ denote the class of measure preserving dynamical
systems having Lebesgue spectrum. Hence, the spectral measure of any nonconstant
$f\in L_2 $ is absolutely continuous with respect to the Lebesgue measure
on $\Delta. $  
\medskip

\noindent{\bf Theorem 2.6}  {\it Let $T\in \mathcal{L} $ and $\ba $ be a
symmetric bounded sequence of complex numbers such that
$ \displaystyle \lim_n \frac1n \sum_{k=0}^{n-1} a_k T^k g(x) $ exists a.e. for every
$g\in L_2 . $  Then $\ba $ is universally good for the eHt in $L_2$ for the class 
$\Cal{L}. $ }

\noindent{\bf Proof.}  First we observe that
$L_2=\Bbb C \oplus \mathcal{H} , $ where
$\mathcal{H}= \{ g\in L_2 : \mu_g =h dz \ \text{for some} \ 0\leq h \leq 1 \} . $
Since constants are trivially good for the eHt, it's enough to prove 
the assertion for $f\in \mathcal{H}.$ 
Again, we will follow the proof of Theorem 2.3;
and hence, it is enough to show that
$ \Vert S_n -S_{-n} \Vert_2 =O(n^{\bt}), $ for some $0<\bt <1, $ where
$S_{\mp j}=\displaystyle \sum_{i=1}^{j} a_i T^{\mp i}f . $
For, using the spectral theorem, we obtain that
$$\Vert S_n -S_{-n} \Vert_2^2 \leq \int_{\Delta}
| \sum_{k=1}^{n} (a_k z^k- a_k \bar{z}^k ) |^2 d\mu_f (z), $$
and hence
$$ \Vert S_n -S_{-n} \Vert_2 \leq
2 [\int_{\Delta} |\sum_{k=1}^{n} a_k z^k|^2 d\mu_f (z)]^{1/2}
+2[\int_{\Delta} |\sum_{k=1}^{n} a_k \bar{z}^k|^2 d\mu_f (z)]^{1/2}. $$
Since $\ba \in l_{\infty}, $ and $T$ has continuous spectrum,
we have,
$$[\int_{\Delta} |\sum_{k=1}^{n} a_k z^k|^2 d\mu_f (z)]^{1/2}
\leq [\int_{\Delta} |\sum_{k=1}^{n} a_k z^k|^2 dz]^{1/2} =
[\Vert f \Vert_2^2 \sum_{k=1}^{n} |a_k |^2 ]^{1/2 }=
O(n^{1/2} \Vert f \Vert_2 ) , $$
and similarly,
$ [\int_{\Delta} |\sum_{k=1}^{n} a_k \bar{z}^k|^2 d\mu_f (z)]^{1/2} =
O( n^{1/2} \Vert f \Vert_2 ) . $
This implies that $ \Vert S_n -S_{-n} \Vert_2 =O(n^{1/2} \Vert f \Vert_2 ) . $ \qed
\medskip

At this point, one might ask if some other subclasses of bounded 
Besicovitch sequences are universally good for the eHt.  
One such candidate is the class of uniform sequences [BK].
However, as the following example shows, they are not the right choice.
Let $Y=\{1,2,3\},\ \ \sigma:Y\to Y $ be a cyclic shift and $\mu $ be the 
uniform $\sigma $-invariant probability measure on $Y$ (i.e., $\mu(\{i\}) =\frac13 $).  Then $(Y,\mu,\sigma)$ is a strictly L-stable system; and hence, 
for each measurable $U\subset Y$ and $y\in Y , $ the sequence $\{a_n\} = \{a_n(y,U)\}$ where
 $$ a_n =\left \{ \aligned & 1 \ \text{if} \ n\geq 0 \ \mathrm{and} \  \sigma ^n(y)\in U, \\ &-1 \ \mathrm{if} \ n<0 \ \mathrm{and} \ a_{-n} = 1, \\ 
&0 \ \ \text{otherwise.} \endaligned \right. $$ Let $U=\{2 \} $ and $y=1, $  then
$$ a_n =\left \{ \aligned & 1 \ \text{if}\ n=3m+1 \ \mathrm{and} \  m\geq 0, \\ & -1 \ \mathrm{if} \ n = 3m-1 \ \mathrm{and} \ m\leq 0, \\
&0 \ \ \text{otherwise.} \endaligned \right. $$
Now, given a measure preserving dynamical system $(X, \Sigma, \nu, \tau), $ let $X$ be divided into 
three sets $A,\ \tau A$ and $\tau^2 A$
with each set having measure $\frac13. $  Define
$$ f(x) =\left \{ \aligned & 0 \ \text{if}\ x\in A \\
& 1 \ \text{if}\ x\in \tau A \\
& -1 \ \text{if}\ x\in \tau^2 A . \endaligned \right. $$
Then, for $x\in A, $ we have,
$$\aligned & f(\tau^{3k} x) =0, \ \ \ f(\tau^{-3k} x) =0 \\
& f(\tau^{3k+1} x) =1, \ \ \ f(\tau^{-3k-1} x) =-1 \\
& f(\tau^{3k+2} x) =-1, \ \ \ f(\tau^{-3k-2} x) =1. \endaligned $$
Therefore, $f\in L_2 $ with $\Vert f \Vert_2 = \sqrt{\frac23} . $
Hence,
$$ \sum_{i=-(3n+1)}^{3n+1} \frac{a_i f(\tau^i x)}{i}
= \sum_{i=1}^{3n+1} \frac{a_i (f(\tau^i x)-f(\tau^{-i} x))}{i} , $$
and for $x\in A $ we have,
$$\sum_{m=1}^{n} \frac{a_{3m+1} (f(\tau^{3m+1} x)-f(\tau^{-3m-1} x))}{3m+1}=
\sum_{m=1}^{n} \frac{1-(-1)}{3m+1}= 2 \sum_{m=1}^{n} \frac{1}{3m+1} . $$
Thus, $\displaystyle \sup_{n\geq 1} 
| \sum_{i=-(3n+1}^{3n+1} \frac{a_i f(\tau^i x)}{i} |=\infty ; $ 
and hence, 
$\displaystyle \lim_n \smm_{i=-n}^{n} \frac{a_i f(\tau^i x)}{i} $
does not exist.
\medskip

\noindent{\bf Remark.}  The system $(X,\Sigma, \nu, \tau) $ 
is not weakly mixing, and the sequence $\{a_n \} $ has 
positive density.
\bigskip

    \bigskip
    
   {\em Fourier coefficients of $L_1$ functions.} In this section we would like to point out that, contrary to the $L_p$ case  $([\c C_2])$, Fourier coefficients of a function $g\in L_1[0,2\pi ]$ need not be universally good for the ergodic Hilbert transform. In fact a much stronger claim holds: there exists a function $g\in L_1[0,2\pi ]$ such that for any dynamical system $(X, \Sigma , \mu , T)$, the Fourier coefficients of $g$ is not good for the ergodic Hilbert transform in $L_p$ for any $p\geq 1$.  
   
   \medskip
   
    We will make use of the fact that Fourier coefficients of the functions from $L_1[0,2\pi]$ may converge to zero arbitrarily slowly. Examples of such functions are well known (a somewhat implicit example can be found in [Z], combining the results in Chapter V, p.183-184 and Chapter III, p.93.). We construct more direct, and somewhat different example (with a not necessarily convex sequence of Fourier coefficients). 
    
    \medskip
    
    {\bf Proposition 2.7.} Let $h :\mathbb{N}\cup \{0\}\rightarrow \mathbb{R}$ be a function such that $\lim _{n\rightarrow \infty }h(n) = 0$. Then there exists a sequence $\{a_n\}_{n\geq 0}$ such that $a_n\geq h (n)$ for every $n\in \mathbb{N}\cup \{0\}$, moreover, for all $x\in [0,2\pi ]$, the series $\frac{1}{2}a_0 + \displaystyle \sum _{n=1}^{\infty }a_ncos(nx) (\star )$  converges to $g(x)$, where $g\in L^1[0,2\pi ]$, and $(\star )$ is the Fourier series of $g$.  
    
    \medskip
    
     For the proof of this proposition we will use the following common notations:
     
     \medskip
     
     {\bf 1.} Given a sequence $\{x_n\}_{n\geq 0}$, we will write $\Delta x_n = x_n - x_{n+1}$, and $\Delta ^2 x_n = \Delta x_n - \Delta x_{n+1}$, for all $n\geq 0$.
     
     \medskip
     
     {\bf 2.} For a function $g\in L^1[0,2\pi ]$, we will write $g^{+} = g\chi _{\{x : g(x)\geq 0\}}$, and $g^{-} = g\chi _{\{x :  g(x)\leq 0\}}$. So $g^{+}$ and $g^{-}$ denote the positive part and the negative part of $g$ respectively. Notice that $g(x) = g^{+}(x) + g^{-}(x), \ \mathrm{a.e.} \  x\in [0,2\pi ]$.
     
     \medskip
     
     {\bf 3.} For all $n\in \mathbb{N}, x\in [0,2\pi ]$ we will write $D_n(x) = \frac{1}{2} + \displaystyle \sum _{k=1}^{n}cos(kx)$ and $F_n(x) = \frac{1}{n+1}\displaystyle \sum _{k=1}^{n}D_n(x)$. 
     
     \medskip
     
     Notice that for all $x\in (0,2\pi )$, $$D_n(x) = \frac{sin(n+\frac{1}{2})x}{2sin(\frac{1}{2}x)} \ \mathrm{and} \ F_n(x) = \frac{2}{n+1}(\frac{sin\frac{1}{2}(n+\frac{1}{2})x}{2sin(\frac{1}{2}x)})^2$$ 
     
     \medskip
     
    {\bf Proof.}  Let $M = 2\mathrm{max}\{h(n) \ : \ n\in \mathbb{N}\cup \{0\}\}$. Since $\displaystyle \lim _{n\rightarrow \infty }h(n) = 0$, we can choose a piecewise linear function $a:[0, \infty )\to (0, \infty )$ satisfying the following conditions
    
    \medskip
    
    (i) $a(n) > h(n), \forall n\in \mathbb{N}$;
    
    \medskip
    
    (ii) $a$ is strictly decreasing;
    
    \medskip
    
    (iii) $a(0) = M$ and $\displaystyle \lim _{x\to \infty }a(x) = 0$;
    
    \medskip
    
    (iv) $a$ is differentiable on $(0,\infty )$ except at countably many points $n_1, n_2, \ldots $, where (assuming $n_0 = 0$), for all $k\in \mathbb{N}$, $n_k$ is an integer and $n_{k} \geq n_{k-1}+3$;
    
    \medskip
    
    (v) there exists $\lambda > 1$ such that for all $k\in \mathbb{N}\cup \{0\}$, $n_{k+1}\leq \lambda (n_{k+1}-n_k)$;
    
    \medskip
    
    (vi) if $s_k$ denotes the slope of $a$ on the interval $(n_{k-1}, n_{k})$ [i.e. $a'(x) = s_k, \forall x\in (n_{k-1}, n_k)$], then for all $k\in \mathbb{N}$, $s_k < s_{k+1}-s_k < -s_k$.

    \medskip
    
  Then we define our sequence $\{a_n\}_{n\geq 0}$ by letting $a_n = a(n), \forall n\in \mathbb{N}\cup \{0\}$.    
    
   \medskip
   
   By definition of the sequence, if $n_{k-1}\leq n \leq n_{k}-2$ then $\Delta ^2a_n = 0$, and if $n = n_{k}-1$ then $|\Delta ^2a_n| = |s_k - s_{k+1}|$.   
    
   \medskip
   
   Then, by the conditions (i)-(vi), we obtain that $\displaystyle \sum _{n=1}^{\infty }(n+1)|\Delta ^2a_n| < \infty  \ (\star _1)$.

   \medskip
   
   Indeed, by conditions (ii), (iii) and (iv), we have $\displaystyle \sum _{i=1}^{\infty }(-s_i)(n_i-n_{i-1}) = M$. Notice that $s_i < 0, \forall i\in \mathbb{N}$. Then, by condition (v), we obtain that $\displaystyle \sum _{i=1}^{\infty }(-s_i)n_i \leq \lambda \displaystyle \sum _{i=1}^{\infty }(-s_i)(n_i-n_{i-1}) < \infty $. This, combined with condition (vi), implies that $\sum _{i=1}^{\infty }|s_i-s_{i+1}|n_i < \infty $; thus we obtain $(\star _1)$. 
    
   \medskip
   
   Let now $s_n(x) = \frac{1}{2}a_0 + \displaystyle \sum _{k=1}^{n}a_kcos(kx)$ for all $n\geq 0, x\in [0,2\pi ]$. Since the sequence $\{a_n\}_{n\geq 0}$ is positive and decreasing, by Abel's summation formula, the limit $\displaystyle \lim _{n\rightarrow \infty }s_n(x)$ exists for all $x\in (0, 2\pi )$. 
   
   \medskip
   
   Also, if $x\in (0,2\pi )$ and $n\geq 1$, by Abel's summation formula, we obtain that $$s_n(x) = \displaystyle \sum _{k=0}^{n-1}\Delta a_kD_k(x) + a_nD_n(x)$$
    Applying Abel's summation formula again, for $n\geq 2$, we get $$s_n(x) = \displaystyle \sum _{k=0}^{n-2}(k+1)\Delta ^2a_kF_k(x) + nF_{n-1}(x)\Delta a_{n-1} + a_nD_n(x)$$
   
   \medskip
   
   Notice that $\displaystyle \lim _{n\rightarrow \infty }a_nD_n(x) = \lim _{n\rightarrow \infty }nF_{n-1}(x)\Delta a_{n-1} = 0$. Hence $g(x) = \displaystyle \lim _{n\rightarrow \infty } \sum _{k=0}^{n-2}(k+1)\Delta ^2a_kF_k(x)$. 
   
   \medskip
   
   Let us now show that $g$ belongs to $L^1[0,2\pi ]$. For all $n\geq 2$, let $$g_n(x) = \sum _{k=0}^{n-2}(k+1)\Delta ^2a_kF_k(x),$$ \ $$S_{+}(n) = \{k \ : \ 0\leq k\leq n, \Delta ^2a_n > 0\}, S_{-}(n) = \{k \ : \ 0\leq k\leq n, \Delta ^2a_k < 0\},$$ \ $$\phi _n(x) = \sum _{k\in S_{+}(n)}(k+1)\Delta ^2a_kF_k(x), \psi _n(x) = \sum _{k\in S_{-}(n)}(k+1)\Delta ^2a_kF_k(x)$$
   
   \medskip
   
   Notice that $F_n(x)\geq 0, \forall x\in [0,2\pi ]$. Then for all $x\in [0,2\pi ]$ we have $g_n(x) = \phi _n(x) + \psi _n(x)$, where  $\phi _n(x)\geq 0$ and $\psi _n(x)\leq 0$. Then $\phi _n(x)\geq g_n^{+}(x)\geq g_n(x), \forall x\in [0,2\pi ]$. By Fatou's Lemma, we have $\displaystyle \int_0^{2\pi }g^{+} = \displaystyle \int_0^{2\pi }\lim _n g_n^{+} \leq \displaystyle \int_0^{2\pi }\liminf _n \phi _n \leq \liminf _n\int_0^{2\pi } \phi _n$. On the other hand, $\int _{0}^{2\pi }F_k(x)dx  = \pi $, for all $k\geq 1$. Then by the condition $(\star _1)$ we obtain that $\displaystyle \liminf _n \int_0^{2\pi } \phi _n < \infty $. Hence $g^{+}\in L^1[0,2\pi ]$. Similarly, we obtain that $g^{-}\in L^1[0,2\pi ]$. Thus $g\in L^1[0,2\pi ]$. 
   
   \medskip
   
  After establishing integrability of $g$ it follows from the claim {\bf 1-8} ([Z], p.184]) that the series $\frac{1}{2}a_0 + \displaystyle \sum _{n=1}^{\infty }a_ncos(nx) (\star )$ is indeed a Fourier series of $g$. $\square $  
         
  \bigskip

  Let now $(X,\Sigma , \mu , T)$ be a dynamical system, and $g\in L^1[0,2\pi ]$ be a function with Fourier coefficients satisfying the following condition: $a_n = 0, \forall n\leq 0$, and $a_n\geq \frac{1}{\log(n)}, \forall n>0$. Let also $\mathbf{a} = \{a_n\}_{n\in \mathbb{Z}}$ and $f(x) = 1, \forall x\in X$. Then $H_{\mathbf{a}}f(x) \geq \displaystyle \sum _{n=1}^{\infty }\frac{1}{n\log(n)} = \infty $, for all $x\in X$. Thus for any dynamical system  $(X,\Sigma , \mu , T)$, the sequence $\{a_n\}$ is not good for eHt in $L_p$ for any $p\geq 1$.  
   
    \bigskip
     
   {\bf Remark.} For convenience of the reader,  we quote the claim {\bf 1-8} from [Z]: {\em If $\{a_n\}_{n\geq 1}$ is a sequence of real numbers decreasing to zero and the function $g(x) = \displaystyle \sum _{n=1}^{\infty }a_ncos(nx)$ is integrable, then the series $\displaystyle \sum _{n=1}^{\infty }a_ncos(nx)$ is a Fourier series of $g$}.
   
\bigskip   

 {\bf Remark.} It is indeed easy to construct a function $a(x)$ with the desired properties (i)-(vi). For all $k\in \mathbb{N}$, let $m_k = \mathrm{min} \{n\in \mathbb{N}\cup \{0\} : h(x) \leq \frac{M}{2^{k+1}}, \forall x \geq n\}$. 
Then let $\{n_k\}_{k\geq 0}$ be a sequence of non-negative integers such that $n_0 = 0, n_k \geq m_k$ and $n_{k+1}\geq 2n_k+3$, for all $k\in \mathbb{N}\cup \{0\}$. We define the function $a(x)$ as follows: we let $a(n_k) = \frac{M}{2^k}, \forall k\in \mathbb{N}\cup \{0\}$, and we affinely extend the function to the interval $[n_k, n_{k+1}]$ for every $k\in \mathbb{N}\cup \{0\}$. Then, by taking $\lambda = 2$, it is clear that all of the conditions (i)-(vi) hold. Notice that the function $a(x)$ constructed in this way will be convex. In general though, conditions (i)-(vi) allow plenty of non-convex functions as well.

 \bigskip

\noindent{\bf 3. Extension to Admissible Processes.}  Given a sequence 
$\ba \in l_{\infty} ,$ let
$$ \Vert \ba \Vert_{\al} := \limsup_{n\geq 1}
\frac{\log^{\al}n}{n^{\al-1}} \sum_{-n}^{n} | a_k | <\infty . $$
Then $\Vert \ \Vert_{\al} $ defines a seminorm on $M_{\al};$
that is $(M_{\al}, \Vert \ \Vert_{\al}) $ is a seminormed subspace of
$l_{\infty}. $  Now, we turn to obtaining some properties of convergence
with respect to $\Vert \ \Vert_{\al} $-seminorm, which will be
instrumental in enlarging the scope of some family of good
modulating sequences.

\noindent{\bf Definition.} A sequence $\ba =\{a_k\}_{k=-\infty }^{\infty} $
of complex numbers is called a {\it Hilbert sequence} if
$\lim_n \smm_{k=-n}^{n} \frac1k a_k $ exists.

\noindent{\bf Remark.}  For any $\lambda\in \Bbb C,\
|\lambda |=1, $ the sequence
$\{\lambda^k \} $ is a Hilbert sequence, and hence, every sequence
induced by a trigonometric polynomial is a (bounded) Hilbert sequence.

\noindent{\bf Proposition 3.1}  {\it a) If $\{ \ba^r \} $ is a
Hilbert sequence for each $r\in \Bbb Z^+ $ and if
$\Vert \ba^r - \ba \Vert_{\al} \to 0 $
as $r \to \infty, $ then $\ba=\{a_k \} $ is a Hilbert sequence.
 
 b) If $\ba^r, \ba \in M_{\al} $ with
$\Vert \ba^r - \ba \Vert_{\al} \to 0 $  and $\ba^r \bb $ is a
bounded Hilbert sequence for all $r\in \Bbb Z^+ , $
for some $\bb \in l_{\infty}, $
then $\ba \bb $ is a bounded Hilbert sequence. }

\noindent{\bf Proof. }  Since
$$\smm_{k=-n}^{n} \frac1k a_k= \smm_{k=-n}^{n}
\frac1k (a_k-a_k^r) +\smm_{k=-n}^{n} \frac1k a_k^r, $$
and since $ \{\smm_{k=-n}^{n} \frac1k a_k^r \}_n $ converges, 
in order to prove (a) it is enough to show that
$\{ \smm_{k=-n}^{n} \frac1k (a_k-a_k^r)\}_n $ converges.
Now, by Abel's summation by parts formula, if
$S^r_{\mp k}=\sum_{i=1}^{k} (a_{\mp i}-a_{\mp i}^r), $ then, for $1\leq m<n, $
$$ \aligned |\smm_{k=-n}^{n} & \frac1k (a_k-a_k^r)
-\smm_{k=-m}^{m} \frac1k (a_k-a_k^r) |
\leq \sum_{m <|k| \leq n} \frac{1}{k(k+1)} |S_k^r | \\
&+|\frac1n S_n^r -\frac{1}{m+1} S_m^r |+|\frac1n S^r_{-n }-\frac{1}{m+1} S^r_{-m }| .
\endaligned $$
By hypothesis, $\Vert \ba^r - \ba \Vert_{\al} \to 0 , $
therefore, given $\ep>0, $ we can pick $N$ large enough such that
whenever $n,r > N, $ we have $\Vert \ba^r - \ba \Vert_{\al} < 1$ and $\displaystyle \sum_{k=n}^{\infty} \frac{1}{k^{3-\al} \log^{\al}k} < \frac {\ep }{2}$.

 Then $|\frac1n S^r_n| \leq
\frac{1}{n}  \smm_{i=-n}^{n} |a_i-a_i^r| \leq
\frac{1}{n} \frac{n^{\al-1}}{\log^{\al}n} \leq \frac{1}{n^{2-\al} \log^{\al}n} , $ hence we have $\displaystyle \lim_n \frac1n S^r_n \to 0.$  Therefore, as $m, n\to \infty, $
$|\frac1n S^r_n-\frac{1}{m+1} S^r_m | \to 0, $ and similarly,
$|\frac1n S^r_{-n}-\frac{1}{m+1} S^r_{-m} | \to 0. $
Hence, for $m,n,r >N $ (by choosing even larger $N $, if necessary),
$$ \sum_{m <|k| \leq n} \frac{1}{k(k+1)} |S_k^r | < \frac{\ep}{2},\ \
|\frac1n S^r_n-\frac{1}{m+1} S^r_m |< \frac{\ep}{4} \ \ \text{and}\ \
|\frac1n S^r_{-n}-\frac{1}{m+1} S^r_{-m} | < \frac{\ep}{4}. $$
Then, it follows that, for $m,n >N, $
$$ |\smm_{k=-n}^{n} \frac1k (a_k-a_k^r)-
\smm_{k=-m}^{m} \frac1k (a_k-a_k^r) |<
\frac{\ep}{2} + \frac{\ep}{4} +\frac{\ep}{4}< \ep, $$
and hence $\{ \smm_{k=-n}^{n} \frac1k (a_k-a_k^r) \}_n $ is Cauchy,
and hence, converges.

Since
$$ \smm_{k=-n}^{n} \frac1k a_k b_k =
\smm_{k=-n}^{n} \frac1k ( a_k-a_k^r) b_k +
\smm_{k=-n}^{n} \frac1k a_k^r  b_k, $$
to prove (b) it is enough to show that
$ |\smm_{k=-n}^{n} \frac1k ( a_k-a_k^r) b_k - \smm_{k=-m}^{m} \frac1k ( a_k-a_k^r) b_k|\to 0 $ as $m, n\to \infty $.
Now, leeting $n > m$, by the inequality
$$|\sum_{m<|k|\leq n} \frac1k ( a_k-a_k^r) b_k |
\leq  \Vert \bb \Vert_{\infty} \sum_{m<|k|\leq n} \frac1k |a_k-a_k^r|, $$
the same method used in part (a) proves the assertion.
\qed

\noindent{\bf Remark.}  Since for any $\lambda\in \Bbb C,\
|\lambda |=1,$ the sequence $\{\lambda^k \} $ is a Hilbert sequence,
for any $\ba \in M_{\al} $ and for
any $ |\lambda |=1 $ the sequence $\{\lambda^k a_k \}$ belongs to $M_{\al} $, and is a Hilbert sequence. 

\medskip

As an application of Proposition 3.1 we will extend the assertion of
Theorem 2.3 to $T$-admissible processes.  Let $(X,\Sigma ,\mu, T)$
be a measure preserving system.  A family of functions
$ \ F=\{ f_{i} \}_{i \in \Bbb Z} \subset L_p(X) , \ 1\leq p \leq \infty, $
is called a {\it T-admissible process on $\Bbb Z $} if
$ T^{\pm{1}} f_{\pm{i}} \leq f_{\pm{(i+1)}} $ for $ i\geq 0. $
When the equality holds, $F$ is called a {\it T-additive process} and 
is necessarily of the form $F=\{ T^i f \}_{i \in \Bbb Z} , $ for some
$f \in L_p(X). $
A process $F=\{ f_{i} \} \subset L_p $ is called
{\it strongly bounded} when
$\sup_{n \in \Bbb Z} \Vert f_n \Vert_p <\infty; $
and it is called {\it symmetric} if
$ T^{2i} f_{-i}=f_i \ \text{for all}\ \ i\in \Bbb Z. $

Given a process $F=\{f_i\} , $ define the Hilbert
transform of $F $ by $\lim_n H_n F(x), $ where
$H_n F(x) =\smm_{i=-n}^{n} \frac1i f_i (x). $
The eHt of a symmetric strongly bounded
T-admissible process $F$ exists a.e. for all $F \subset L_1 $ [\c C$_1$].
There, it is also shown that if $F=\{f_n\}\subset L_p $ is a positive
symmetric strongly bounded $T$-admissible process, then
there exists a monotone increasing sequence
$\{v_r \} \in L_p^+ $ and $v_r \uparrow \dl \in L_p$ such that
$f_n=T^n v_{|n|} $ for all $n\in \Bbb Z , $
$f_n \leq T^n \dl $ for all $n\in \Bbb Z , $ and
$ \Vert \dl \Vert_p =\sup_{n\in \Bbb Z} \Vert f_n \Vert_p .  $

For $r \geq 1 , $ define
$ g_i^r(x)= f_i (x) $ for $ 0\le |i| \le r $ and
$$ g_i^r(x) \ =
\left\{ \aligned & T^{i-r}f_r(x) \ \ \text{for} \ \ i > r \\
& T^{-i+r}f_{-r}(x) \ \ \text{for} \ \ -i > r . \endaligned \right.  $$
Thus, $g_i^r(x) \leq f_i (x) $ for every $i\in \Bbb Z $ and
for each $r\geq 1, $  and,
$$ 0\leq f_i(x)-g_i^r(x) \leq T^i (\dl - v_r )(x) \ \ \text{if}\ |i|>r,  \ \
\text{and} \ 0 \ \text{if}\ |i|\leq r. $$
Observe that, $\Vert \dl - v_r \Vert_p \downarrow 0  $ as $r\to \infty. $
Furthermore, ignoring the first $r$ terms, the process $\{ g_k^r \}_k $ is 
$T$-additive.  It follows that, if $\ba \in \Bbb C \oplus A_{\al } , $ 
then, for each $r \geq 1, $
$\displaystyle \lim_n \smm_{-n}^n \frac{a_i g_i^r}{i} $ exists a.e. by Theorem 2.3.
Therefore, for each $r\geq 1, $ for a.e. $x\in X, $ the sequence
$\{a_i g_i^r (x) \} $ is a Hilbert sequence.
\medskip

\noindent{\bf Theorem 3.2}  {\it Let $F \subset L_2 $ be a symmetric,
strongly bounded $T$-admissible process.
If $\ba \in  M_{\al }, $ then  
$$ \lim_n \smm_{-n}^n \frac{a_i f_i(x)}{i} \ \ \text{exists a.e.} $$ }

\noindent{\bf Proof.} Let $\bu^r $ and $\bu $ be defined by
$\bu^r =\{ a_i g_i^r \} $ and $ \bu =\{a_i f_i \} , \ r\geq 1 , $
respectively.  By the assumptions, for each $r, $ the sequence $\bu^r $
is a Hilbert sequence a.e. $x \in X. $  Since,
$$ 0\leq \frac{\log^{\al}n}{n^{\al-1}} \sum_{i=-n}^{n} |a_i g_i^r - a_i f_i |
\leq \frac{\log^{\al}n}{n^{\al-1}} \sum_{i=-n}^{n} |a_i| T^i (\dl - v_r) , $$
it follows that,
$$ 0\leq \int [\frac{\log^{\al}n}{n^{\al-1}} \sum_{i=-n}^{n} |a_i| T^i (\dl - v_r) ] d\mu
\leq \frac{\log^{\al}n}{n^{\al-1}} \sum_{i=-n}^{n} |a_i| \Vert \dl - v_r \Vert_2
\leq C_{\ba} \Vert \dl - v_r \Vert_2 \to 0. $$
Hence, by Proposition 3.1 (a), it follows that
$\bu=\{a_i f_i (x) \} $ is Hilbert sequence for a.e., which
proves that $ \displaystyle \lim_n \smm_{-n}^n \frac{a_i f_i(x)}{i} $ exists a.e. \qed
\medskip

\noindent{\bf Corollary 3.3}  {\it Let $(X, \Sigma, \mu, T)$ be
a measure preserving system and
$\ba \in MB_{\al }$ be a two-sided sequence.  Then
$$ \lim_n \smm_{-n}^n \frac{a_i f_i(x)}{i} \ \ \text{exists a.e.} $$
for any symmetric, strongly bounded $T$-admissible process
$F=\{f_k\}  \subset L_2(X) . $  If $\ba \in MB_{\al } $ is a one-sided sequence, 
then the assertion holds if $(X, \Sigma, \mu, T)$ is a weakly mixing system.}

\vspace{2cm}

\centerline{\bf REFERENCES}

\def\ref{\leftskip 20pt \parindent-20pt\parskip 4pt}

\ref [A] I. Assani, {\it A Wiener-Wintner property for the helical transform,}
Ergod. Th. \& Dynam. Syst., {\bf 12}, 185-194, 1992.

\ref [BL], A. Bellow and V. Losert, {\it The weighted pointwise ergodic theorem 
and the inividual ergodic theorem along subsequences}, Trans. AMS, {\bf 288},
307-345, 1985.

\ref [BK] A. Brunel and M. Keane, {\it Ergodic theorems for operator sequences,}
Zeit. Wahr., {\bf 12}, 231-240, 1969.

\ref [C] M. Cotlar, {\it A unified theory of of Hilbert transforms and ergodic
theorem}, Rev. mat. Cuyana, {\bf 1}, 105-167, 1955.

\ref [CP] J. Campbell and K. Petersen, {\it The spectral measure and Hilbert 
transform of a measure-preserving transformation}, Trans. AMS, {\bf 313},
121-129, 1989.

\ref [\c C$_1$] D.\c C\"omez, {\it Existence of discrete ergodic
singular transforms for admissible processes}, Colloq. Math., {\bf 112},
335-343, 2008.

\ref [\c C$_2$] D.\c C\"omez, {\it The modulated ergodic Hilbert transform},
Discrete \& Cont. Dynam. Syst., Series S, {\bf 2}, 325-336, 2009. 

\ref [J] R. Jajte, {\it On the existence of the ergodic Hilbert transform},
The Annals of Prob., {\bf 15}, 831-835, 1987

\ref [K] J.P. Kahane, {\it Sur les coefficients de Fourier-Bohr}, Studia Math.,
{\bf 21}, 103-106, 1961.

\ref [LT] M. Lacey and E. Terwilleger, {\it A Wiener-Wintner theorem for the Hilbert
transform}, Ark. Math., {\bf 46}, 315-336, 2008.

\ref [P] K.Petersen, {\it Another proof of the existence of the ergodic
Hilbert transform,} Proc. AMS, {\bf 88}, 39-43, 1983.

\ref [S] R. Sato, {\it A remark on the ergodic Hilbert transform}, 
Math. J. Okayama Univ., {\bf 28}, 159-163, 1986.

\ref [Z] A. Zygmund, {\it Trigonometric Series.} Cambridge University Press 1959.

\bigskip

\noindent Azer Akhmedov 

\noindent Department of Mathematics, 
North Dakota State University, 
Department \# 2750, 
PO Box 6050,
Fargo, ND 58108-6050, 
USA. 

\noindent {\it E-mail address:} azer.akhmedov@ndsu.edu

\bigskip

\noindent Do\u gan \c C\"omez 

\noindent Department of Mathematics, 
North Dakota State University, 
Department \# 2750, 
PO Box 6050,
Fargo, ND 58108-6050, 
USA. 

\noindent {\it E-mail address:} Dogan.Comez@ndsu.edu

\end{document}